\documentclass{article}
\usepackage{indentfirst}

\usepackage{graphicx}
\usepackage{epstopdf}
\usepackage{float}
\usepackage{subfig}
\usepackage{amssymb}
\usepackage{stmaryrd}
\usepackage{amsfonts}
\usepackage{amsmath}
\usepackage{latexsym}
\usepackage{color}
\usepackage{CJK}
\usepackage{cite}
\usepackage{mathrsfs}
\usepackage{algorithm} 
\usepackage{algorithmic} 
\usepackage{multirow} 
\usepackage{amsmath}
\usepackage{xcolor}
\usepackage{graphicx}
\usepackage{mathptmx}
\usepackage{indentfirst}
\usepackage{booktabs}
\usepackage{float}
\usepackage{subfig}
\usepackage{amssymb}
\usepackage{stmaryrd}
\usepackage{dsfont}
\usepackage{amssymb}
\usepackage{amsfonts}
\usepackage{amsmath}
\usepackage{latexsym}
\usepackage{color}
\usepackage{cite}
\usepackage{mathrsfs}
\usepackage{epstopdf}

\begin{document}

\title{{\bf Stability Analysis of Short Memory Fractional Differential Equations }}
\author{Xudong Hai$^1$, Guojian Ren$^1$, Yongguang Yu$^{1,}$\thanks{Corresponding author
E-mail: ygyu@bjtu.edu.cn}, Lipo Mo$^2$, Conghui Xu$^1$}
\date{\small $^{1.}$ Department of Mathematics, Beijing Jiaotong University,
Beijing, 100044, P.R.China\\
\small $^{2.}$ Department of Mathematics, Beijing Technology and Business University,
Beijing, 100048, P.R.China}

\maketitle

\begin{center}
\begin{minipage}{13cm}\vspace{0.1cm}
{\bf Abstract:} In this paper, a fractional derivative with short-term memory properties is defined, which can be viewed as an extension of Caputo fractional derivative. Then, some properties of the short memory fractional derivative are discussed. Also, a comparison theorem for a class of short memory fractional systems is shown, via which some relationship between short memory fractional systems and Caputo fractional systems can be established. By applying the comparison theorem and Lyapunov direct method, some sufficient criteria are obtained, which can ensure the asymptotic stability of some short memory fractional equations. Moreover, a special result is presented, by which the stability of some special systems can be judged directly. Finally, three examples are provided to demonstrate the effectiveness of the main results.



\end{minipage}
\end{center}
\vspace{0.5cm}

\section{Introduction}
Fractional calculus \cite{1} is a classic mathematical concept, which can be viewed as an generalization from integer order differentiation and integration to arbitrary non-integer order, and it has become very popular in recent years due to its potential application value in many fields of science and engineering \cite{2}. For example, it has been found that many real-world systems and processes can be modelled via using fractional derivatives such as diffusion and wave propagation \cite{3,4}, dielectric relaxation phenomena in polymeric materials \cite{5}, boundary layer effects in ducts \cite{6} and so on. And more and more efforts are being made to enrich the theories and applications of fractional calculus.

Different from classical integer order systems, fractional systems can provide an excellent instrument for the description of memory and hereditary properties of various materials and processes. The main reason for this difference is that the definition of fractional derivative is different from that of integer order derivative. For example, from the definition of Caputo fractional derivative \cite{2}, it is easy to observe that the fractional derivative of a function at any given time depends on all the historical information before the moment, but the integer order derivative does not. This special property of fractional derivatives makes it possible to use fractional systems to model some real-world systems with memory properties, and some efforts have been devoted \cite{7,8}. However, the common fractional derivatives, like Caputo fractional derivative and Riemann-Liouville fractional derivative, are operators with long-term memory properties, because these differential operators rely on all of the historical information rather than some of the historical information. It will be difficult to use these fractional derivatives to describe some processes with short-term memory properties, like visual process \cite{9}. Thus it is necessary to design and study some special differential operators with short memory properties. Fortunately, some attempts have been made \cite{10,11}.

Similar to integer order differential systems, the investigation of stability is always an important task for fractional systems. Over the past few decades, a lot of results about the stability analysis of fractional systems have been obtained. And some examples can be cited. The generalized Mittag-Leffler stability of fractional nonlinear dynamic systems was studied in \cite{12}. The stability of fractional neural networks was fully investigated through an energy-like function analysis in \cite{13}. A comparison theorem for a class of fractional systems with time delays was proposed in \cite{14}.  A new property for Caputo fractional derivative was presented in \cite{15}. Some sufficient conditions were established to ensure the existence and uniqueness of the nontrivial solution for some fractional equations in \cite{16}. The LMI-based stability of non-autonomous fractional systems with multiple time delays was investigated in \cite{17}. Of course, there are still a lot of investigations on this topic being carried out and it is conceivable that many excellent results will be obtained in the future. However, to the best of our knowledge, there is no known result concerning the stability analysis for fractional systems with short-term memory properties.

Motivated by the above discussion, the main contributions of this paper can be summarized as follows. Firstly, a new fractional derivative with short-term memory properties is defined, which can be viewed as an extension of Caputo fractional derivative, and some of its properties are proved. Secondly, a comparison theorem for short memory fractional systems is proposed, via which some relationship between short memory fractional systems and Caputo fractional systems are established. Then some asymptotic stability conditions of short memory fractional equations are obtained. In addition, a special theorem is proposed, according to which the stability of some special differential equations can be judged directly.

The rest of this paper is organized as follows. Section $2$ gives some basic concepts about fractional calculus and several conclusions about the stability analysis of fractional systems. In section $3$, the definition of a short memory fractional derivative is proposed and some stability criteria of short memory fractional systems are derived. Next, three numerical examples are presented to illustrate the correctness of the main results in section $4$. Finally, the conclusion of this paper is given in section $5$.

$Notations$: Let $R^+=[0,+\infty)$, $R=(-\infty,+\infty)$, $R^n$ be the $n$-dimensional Euclidean space and $R^{n\times m}$ be the space of $n\times m$ real matrices. $A^T$ represents the transpose of matrix $A$. $\parallel \cdot \parallel$ denotes the Euclidean norm of a vector, for example, $\parallel \beta\parallel=\sqrt{\sum_{i=1}^n\beta_i^2}$ where $\beta=(\beta_1,\cdots,\beta_n)^T$.

\section{Preliminaries}
In this section, two common definitions of fractional calculus are introduced and several necessary lemmas are presented.

\textbf{Definition 1 \cite{2} }: Caputo fractional derivative with order $\alpha$ for a function $x: R^+\rightarrow R$ is defined as
\[{^C_{{t_0}}D_t^\alpha }x(t)=\frac{1}{\Gamma (m-\alpha)}\int_{{t_0}}^t(t-\tau)^{m-\alpha-1}\frac{d^mx(\tau)}{d\tau^m}d\tau,\]
where $0\leq m-1< \alpha < m$, $m\in Z_+$ and $\Gamma(\cdot)$ is Gamma function.

\textbf{Definition 2 \cite{2} }: Riemann-Liouville fractional integral with order $\alpha$ for a function $f:R^+\rightarrow R $ is defined by
\[{^R_{{t_0}}I_t^\alpha }f\left( t \right) = \frac{1}{{\Gamma \left( {\alpha} \right)}}\int_{{t_0}}^t (t-\tau)^{\alpha-1}f(\tau)d\tau,\]
where $\alpha > 0$ and $\Gamma(\cdot)$ is Gamma function.

\textbf{Lemma 1 \cite{2} }: For any constants $k_1$ and $k_2$ , the linearity of Caputo fractional derivative is described by
\[{^C_{{t_0}}D_t^\alpha }(k_1 g(t)+k_2h(t))=k_1{^C_{{t_0}}D_t^\alpha }g(t)+k_2{^C_{{t_0}}D_t^\alpha }h(t).\]

\textbf{Lemma 2 \cite{2} }: For any constants $k_1$ and $k_2$ , the linearity of Riemann-Liouville fractional integral is described by
\[{^R_{{t_0}}I_t^\alpha }(k_1 g(t)+k_2h(t))=k_1{^R_{{t_0}}I_t^\alpha }g(t)+k_2{^R_{{t_0}}I_t^\alpha }h(t).\]

\textbf{Lemma 3 \cite{2} }: If $0<\alpha<1$ and $x(t):R\rightarrow R$ is a differentiable function, then the following equality holds
\[{^R_{t_0}I^\alpha_{t}}{^C_{t_0}D^\alpha_t}x(t)=x(t)-x(t_0),\]
where $t\geq t_0$.

\textbf{Lemma 4 \cite{14} }: Consider the following fractional differential inequality with time delay:
\[
  \left\{
\begin{aligned}
&{^C_{t_0}D^\alpha_t}x(t)\leq -ax(t)+bx(t-q),~~0<\alpha\leq 1,\\
&x(t)=h(t),~~t\in[t_0-q,t_0],
  \end{aligned}
\right.
\]
and the linear fractional differential systems with time delay
\[\left\{
\begin{aligned}
&{^C_{t_0}D^\alpha_t}y(t)= -ay(t)+by(t-q),~~0<\alpha\leq 1,\\
&y(t)=h(t),~~t\in[t_0-q,t_0],
  \end{aligned}
\right.\]
where $x(t)\in R$ and $y(t)\in R$ are continuous and nonnegative in $[t_0,+\infty)$, $q$ is a positive constant and $h(t)\geq 0$ for $t\in[t_0-q,t_0]$. If $a>0$ and $b>0$, then $x(t)\leq y(t)$ for $t\in[t_0,+\infty)$.

\textbf{Lemma 5 \cite{14} }: Consider the following linear fractional system with time delay:
\begin{equation}\left\{
\begin{aligned}
&{^C_{t_0}D^\alpha_t}x(t)=-ax(t)+bx(t-q),~~0<\alpha <1,\\
&x(t)=h(t),~~t\in[t_0-q,t_0],
\end{aligned}\right.
\end{equation}
where $x(t)\in R$, $q$ is a positive constant and $h(t)$ is bounded in $[t_0-q,t_0]$. If $a>0$, $b>0$ and $a>b$, then the zero equilibrium point of the system $(1)$ is Lyapunov asymptotically stable.

Consider the Caputo fractional non-autonomous system
\begin{equation}
\begin{aligned}
{^C_{t_0}D^\alpha_t}x(t)=f(x,t)
\end{aligned}
\end{equation}
with initial condition $x(t_0)$, where $0<\alpha<1$, $x(t)\in R^n$, $f:[t_0,+\infty)\times \Omega\rightarrow R^n$ is continuous in $t$ and $x$ on $[t_0,+\infty)\times \Omega$, and $\Omega\subset R^n$ is a domain that contains the origin $x=0$.

\textbf{Definition 3 \cite{12} }: The constant $x_0$ is an equilibrium point of the Caputo fractional system $(2)$, if and only if $f(t,x_0)=0$.

\textbf{Definition 4 \cite{18} }: Let $\phi(\gamma):R^+\rightarrow R^+$ be a continuous function, where $R^+\in \{\gamma\mid \gamma\geq 0\}$. $\phi(\gamma)$ is said to belong to class-$K$ if it is strictly increasing and $\phi(0)=0$.

\textbf{Lemma 6 \cite{12} }: Let $x=0$ be an equilibrium point for the non-autonomous fractional system $(2)$. Assume that there exists
a Lyapunov function $V(t,x(t))$ and class-$K$ functions $\phi_i$ $(i=1,2,3)$ satisfying
\[\phi_1(\parallel x\parallel)\leq V(t,x)\leq \phi_2(\parallel x\parallel)\]
and
\[{^C_{t_0}D^\alpha_t}V(t,x(t))\leq-\phi_3(\parallel x\parallel)\]
where $0<\alpha<1$. Then the system $(2)$ is Lyapunov asymptotically stable.

\textbf{Lemma 7 \cite{15} }: Let $v(t):R\rightarrow R^n$ be a differentiable vector function. Then, for any time instant $t\geq t_0$, the following inequality holds
\begin{equation}
\begin{aligned}
{^C_{{t_0}}D_t^\alpha}(v^T(t)Pv(t))\leq 2v^T(t)(P{^C_{{t_0}}D_t^\alpha}v(t)),
\end{aligned}
\end{equation}
where $P\in R^{n\times n}$ is a positive definite matrix and $0<\alpha <1$.

\section{Main results}
In this section, a special fractional derivative with short-term memory properties is defined and some of its properties are discussed. Then the stability analysis of short memory fractional differential equations is given and some useful conclusions are obtained.

\subsection{Fractional derivative with short-term memory properties}
The definition of the fractional derivative with short-term memory properties is given as follows.

\textbf{Definition 5 }: The short memory fractional derivative of a continuously differentiable function $f(t):R^+\rightarrow R$ is defined as
\begin{equation}
\begin{aligned}
{_{s_\omega(t)}\tilde{D}^\alpha_t}f(t)=\frac{1}{\Gamma(1-\alpha)}\int_{s_\omega(t)}^t(t-\tau)^{-\alpha}\frac{df(\tau)}{d\tau}d\tau,
\end{aligned}
\end{equation}
where $0<\alpha<1$ and $\Gamma(\cdot)$ is Gamma function. And the function $s_\omega(t)$ is described as
\begin{equation}\begin{split}\ s_\omega(t)=\begin{cases}
t_0,&\text{$t_0\leq t\leq t_0+\omega$},\\
t-\omega,
&\text{$t>t_0+\omega$},
\end{cases}
\end{split}
\end{equation}
where $\omega$ is a positive constant.

\textbf{Remark 1 }: Obviously, when $t_0\leq t\leq t_0+\omega$,
\[{_{s_\omega(t)}\tilde{D}^\alpha_t}f(t)={^C_{t_0}D^\alpha_t}f(t).\] And in particular, when $\omega\rightarrow+\infty$,
\[{_{s_\omega(t)}\tilde{D}^\alpha_t}f(t)={^C_{t_0}D^\alpha_t}f(t),\]
which means that the short memory derivative $(4)$ can be seen as an extension of Caputo fractional derivative. Meanwhile, it is easy to get a conclusion that ${_{s_\omega(t)}\tilde{D}^\alpha_t}f(t)$ is a continuous function because $s_\omega(t)$ is continuous for $t\geq t_0$.

\textbf{Remark 2 }: According to $Definition~1$, it is easy to find that Caputo fractional derivative of the function $x$ at time $t$ depends on the information of $x$ on the interval $[t_0,t]$. This means that Caputo fractional derivative is an operator with long-term memory properties. And according to $Definition~5$, we can find that the short memory derivative of the function $f$ at time $t$ relies on the information of $f$ on the interval $[s_\omega(t),t]$. This means that the operator ${_{s_\omega(t)}\tilde{D}^\alpha_t}$ also has memory properties. But compared with Caputo fractional derivative, ${_{s_\omega(t)}\tilde{D}^\alpha_t}$ is a special operator with short-term memory properties.

Then, two conclusions about the short memory derivative $(4)$ are derived.

\textbf{Theorem 1 }: The linearity of the operator ${_{s_\omega(t)}\tilde{D}^\alpha_t}$ can be described as follows.

Consider continuously differentiable functions $f(t)$ and $g(t)$. For any constants $k_1$ and $k_2$, the following equality holds,
\begin{equation}
\begin{aligned}
{_{s_\omega(t)}\tilde{D}^\alpha_t}[k_1f(t)+k_2g(t)]=k_1{_{s_\omega(t)}\tilde{D}^\alpha_t}f(t)+k_2{_{s_\omega(t)}\tilde{D}^\alpha_t}g(t),
\end{aligned}
\end{equation}
where $0<\alpha<1$, $s_\omega(t)$ has been described in the equality $(5)$ and $\omega$ is a positive constant.

\textbf{Proof }: According to $Definition~5$, the following equality can be obtained,
\begin{equation}
\begin{aligned}
{_{s_\omega(t)}\tilde{D}^\alpha_t}[k_1f(t)+k_2g(t)]&=\frac{1}{\Gamma(1-\alpha)}\int_{s_\omega(t)}^t(t-\tau)^{-\alpha}\frac{d(k_1f(\tau)+k_2g(\tau))}{d\tau}d\tau\\
&=k_1\frac{1}{\Gamma(1-\alpha)}\int_{s_\omega(t)}^t(t-\tau)^{-\alpha}\frac{df(\tau)}{d\tau}d\tau+k_2\frac{1}{\Gamma(1-\alpha)}\int_{s_\omega(t)}^t(t-\tau)^{-\alpha}\frac{dg(\tau)}{d\tau}d\tau\\
&=k_1{_{s_\omega(t)}\tilde{D}^\alpha_t}f(t)+k_2{_{s_\omega(t)}\tilde{D}^\alpha_t}g(t).
\end{aligned}
\end{equation}

This completes the proof.

\textbf{Theorem 2 }: Let $x(t)\in R$ be a continuously differentiable function. For any time instant $t\geq t_0$, the following inequality holds,
\begin{equation}
\begin{aligned}
\frac{1}{2}{_{s_\omega(t)}\tilde{D}^\alpha_t}x^2(t)\leq x(t){_{s_\omega(t)}\tilde{D}^\alpha_t}x(t),
\end{aligned}
\end{equation}
where $0<\alpha<1$, $s_\omega(t)$ is described in the equality $(5)$ and $\omega$ is a positive constant.

\textbf{Proof }: When $t_0\leq t\leq t_0+\omega$,
\begin{equation}
\begin{aligned}
\frac{1}{2}{_{s_\omega(t)}\tilde{D}^\alpha_t}x^2(t)&=\frac{1}{2}{^C_{t_0}D^\alpha_t}x^2(t)\\
&\leq x(t){^C_{t_0}D^\alpha_t}x(t)\\
&= x(t){_{s_\omega(t)}\tilde{D}^\alpha_t}x(t),
\end{aligned}
\end{equation}
where $Lemma~7$ is applied above.

When $t> t_0+\omega$, let $y(\tau)=x(t)-x(\tau)$, which implies that $\frac{dy(\tau)}{d\tau}=-\frac{dx(\tau)}{d\tau}$. Then the following equality can be derived,
\begin{equation}
\begin{aligned}
\frac{1}{\Gamma(1-\alpha)}\int_{t-\omega}^t\frac{y(\tau)\frac{dy(\tau)}{d\tau}}{(t-\tau)^\alpha}d\tau =\frac{1}{2\Gamma(1-\alpha)}\{\frac{y^2(\tau)}{(t-\tau)^\alpha}\mid_{\tau=t}-\frac{y^2(t-\omega)}{\omega^\alpha}-\alpha\int_{t-\omega}^t\frac{y^2(\tau)}{(t-\tau)^{\alpha+1}}d\tau\}.
\end{aligned}
\end{equation}

By using L'Hopital rule, we have
\begin{equation}
\begin{aligned}
\frac{y^2(\tau)}{(t-\tau)^\alpha}\mid_{\tau=t}&=\lim\limits_{\tau\rightarrow t-}\frac{(x(t)-x(\tau))^2}{(t-\tau)^\alpha}\\
&=\lim\limits_{\tau\rightarrow t-}\frac{-2(x(t)-x(\tau))\frac{dx(\tau)}{d\tau}}{-\alpha(t-\tau)^{\alpha-1}}\\
&=\lim\limits_{\tau\rightarrow t-}\frac{2}{\alpha}(x(t)\frac{dx(\tau)}{d\tau}-x(\tau)\frac{dx(\tau)}{d\tau})(t-\tau)^{1-\alpha}\\
&=0.
\end{aligned}
\end{equation}

It follows from $(10)$ and $(11)$ that
\begin{equation}
\begin{aligned}
\frac{1}{\Gamma(1-\alpha)}\int_{t-\omega}^t\frac{y(\tau)\frac{dy(\tau)}{d\tau}}{(t-\tau)^\alpha}d\tau\leq 0.
\end{aligned}
\end{equation}

From $Theorem~1$ and the inequality $(12)$,
\begin{equation}
\begin{aligned}
\frac{1}{2}{_{s_\omega(t)}\tilde{D}^\alpha_t}x^2(t)- x(t){_{s_\omega(t)}\tilde{D}^\alpha_t}x(t)&=\frac{1}{\Gamma(1-\alpha)}\int_{t-\omega}^t\frac{x(\tau)\frac{dx(\tau)}{d\tau}-x(t)\frac{dx(\tau)}{d\tau}}{(t-\tau)^\alpha}d\tau\\
&=\frac{1}{\Gamma(1-\alpha)}\int_{t-\omega}^t\frac{-(x(t)-x(\tau))\frac{dx(\tau)}{d\tau}}{(t-\tau)^\alpha}d\tau\\
&=\frac{1}{\Gamma(1-\alpha)}\int_{t-\omega}^t\frac{y(\tau)\frac{dy(\tau)}{d\tau}}{(t-\tau)^\alpha}d\tau\\
&\leq 0.
\end{aligned}
\end{equation}

Then, the following inequality holds for $t\geq t_0$,
\begin{equation}
\begin{aligned}
\frac{1}{2}{_{s_\omega(t)}\tilde{D}^\alpha_t}x^2(t)\leq x(t){_{s_\omega(t)}\tilde{D}^\alpha_t}x(t).
\end{aligned}
\end{equation}

This completes the proof.

\textbf{Remark 3 }: Similar to the proof of $Theorem~2$, the following relation can be derived directly,
\begin{equation}
\begin{aligned}
\frac{1}{2}{_{s_\omega(t)}\tilde{D}^\alpha_t}x^T(t)x(t)\leq x^T(t){_{s_\omega(t)}\tilde{D}^\alpha_t}x(t),
\end{aligned}
\end{equation}
where $0<\alpha<1$ and $x(t)\in R^n$.

\subsection{Stability analysis of short memory fractional differential equations}
An important conclusion is given as follows, which will be used in the derivation of several subsequent theoretical results. And this conclusion can be called a comparison theorem for short memory fractional differential equations.

\textbf{Theorem 3 }: Consider the following short memory fractional differential inequality
\begin{equation}
  \left\{
\begin{aligned}
&{_{s_\omega(t)}\tilde{D}^\alpha_t}x(t)\leq -ax(t),~~0<\alpha<1,\\
&x(t_0)=x_0,
  \end{aligned}
\right.
\end{equation}
and the Caputo fractional differential equation with time delay
\begin{equation}
  \left\{
\begin{aligned}
&{^C_{t_0}D^\alpha_t}y(t)= -ay(t)+\frac{1}{\omega^\alpha\Gamma(1-\alpha)}y(t-\omega),~~0<\alpha<1,\\
&y(s)=x_0,~~s\in[t_0-\omega,t_0],
  \end{aligned}
\right.
\end{equation}
where $x(t)$ and $y(t)$ are continuous and nonnegative in $[t_0,+\infty)$, $y(s)=x_0\geq 0$ for $s\in[t_0-\omega,t_0]$, $s_\omega(t)$ is described in the equality $(5)$, $\omega$ is a positive constant and $\Gamma(\cdot)$ is Gamma function.

If $a>0$, then
\begin{equation}
\begin{aligned}
x(t)\leq y(t),~~t\in[t_0,+\infty).
\end{aligned}
\end{equation}

\textbf{Proof }: When $t_0\leq t\leq t_0+\omega$, the following inequality holds,
\begin{equation}
\begin{aligned}
{^C_{t_0}D^\alpha_t}x(t)={_{s_\omega(t)}\tilde{D}^\alpha_t}x(t)&\leq -ax(t)\\
&\leq-ax(t)+\frac{1}{\omega^\alpha\Gamma(1-\alpha)}x(t_0).
\end{aligned}
\end{equation}

According to $Definition~1$ and $Definition~5$, the following inequality holds when $t>t_0+\omega$,
\begin{equation}
\begin{aligned}
{^C_{t_0}D^\alpha_t}x(t)&={_{s_\omega(t)}\tilde{D}^\alpha_t}x(t)+\frac{1}{\Gamma(1-\alpha)}\int_{t_0}^{s_\omega(t)}(t-\tau)^{-\alpha}\frac{dx(\tau)}{d\tau}d\tau\\
&\leq-ax(t)+\frac{1}{\Gamma(1-\alpha)}\int_{t_0}^{t-\omega}(t-\tau)^{-\alpha}\frac{dx(\tau)}{d\tau}d\tau\\
&=-ax(t)+\frac{1}{\Gamma(1-\alpha)}\{\frac{x(\tau)}{(t-\tau)^\alpha}\mid_{t_0}^{t-\omega}-\int_{t_0}^{t-\omega}x(\tau)d(\frac{1}{(t-\tau)^\alpha})\}\\
&=-ax(t)+\frac{1}{\Gamma(1-\alpha)}\{\frac{x(t-\omega)}{\omega^\alpha}-\frac{x(t_0)}{(t-t_0)^\alpha}-\alpha\int_{t_0}^{t-\omega}(t-\tau)^{-\alpha-1}x(\tau)d\tau\}\\
&\leq -ax(t)+\frac{1}{\omega^\alpha\Gamma(1-\alpha)}x(t-\omega).
\end{aligned}
\end{equation}

From inequalities $(19)$ and $(20)$, the following fractional differential inequality is obtained,
\begin{equation}
  \left\{
\begin{aligned}
&{^C_{t_0}D^\alpha_t}x(t)\leq -ax(t)+\frac{1}{\omega^\alpha\Gamma(1-\alpha)}x(t-\omega),\\
&x(s)=x_0,~~s\in[t_0-\omega,t_0].
  \end{aligned}
\right.
\end{equation}

Consider the following linear fractional differential system with time delay,
\begin{equation}
  \left\{
\begin{aligned}
&{^C_{t_0}D^\alpha_t}y(t)= -ay(t)+\frac{1}{\omega^\alpha\Gamma(1-\alpha)}y(t-\omega),\\
&y(s)=x_0,~~s\in[t_0-\omega,t_0].
  \end{aligned}
\right.
\end{equation}

According to $Lemma~4$, it follows from $(21)$ and $(22)$ that
\begin{equation}
\begin{aligned}
x(t)\leq y(t), ~~t\geq t_0.
\end{aligned}
\end{equation}

This completes the proof.

\textbf{Remark 4 }: $Theorem~3$ is similar to the fractional comparison theorem proposed in \cite{12,14}, which may be useful in deriving some theoretical results.

Consider the short memory fractional differential system
\begin{equation}
  \left\{
\begin{aligned}
&{_{s_\omega(t)}\tilde{D}^\alpha_t}x(t)= f(x,t),~~0<\alpha<1,\\
&x(t_0)=x_0,
  \end{aligned}
\right.
\end{equation}
where $x(t)\in R^n$, $f:[t_0,+\infty)\times R^n\rightarrow R^n$ is continuous in $t$ and $x$ on $[t_0,+\infty)\times R^n$, $s_\omega(t)$ is described in the equality $(5)$ and $\omega$ is a positive constant.

\textbf{Remark 5 }: When $t_0\leq t\leq t_0+\omega$, the following equation holds,
\[
{_{s_\omega(t)}\tilde{D}^\alpha_t}x(t)= {^C_{t_0}D^\alpha_t}x(t)=f(x,t).
\]
And when $t>t_0+\omega$, the following equality can be derived,
\[\begin{aligned}
{^C_{t_0}D^\alpha_t}x(t)&={^C_{t-\omega}D^\alpha_t}x(t)+\frac{1}{\Gamma(1-\alpha)}\int_{t_0}^{t-\omega}(t-\tau)^{-\alpha}\frac{dx(\tau)}{d\tau}d\tau\\
&={_{s_\omega(t)}\tilde{D}^\alpha_t}x(t)+\frac{1}{\Gamma(1-\alpha)}\int_{t_0}^{t-\omega}(t-\tau)^{-\alpha}\frac{dx(\tau)}{d\tau}d\tau\\
&=f(x,t)+\frac{1}{\Gamma(1-\alpha)}\{\frac{x(t-\omega)}{\omega^\alpha}-\frac{x(t_0)}{(t-t_0)^\alpha}-\alpha\int_{t_0}^{t-\omega}(t-\tau)^{-\alpha-1}x(\tau)d\tau\}.
\end{aligned}\]
Obviously, the solution $x(t)$ of the equation $(24)$ also satisfies the following differential system,
\begin{equation}\begin{split}\
{^C_{t_0}D^\alpha_t}x(t)=\begin{cases}
f(x,t)&\text{$t_0\leq t\leq t_0+\omega$},\\
f(x,t)+\frac{1}{\Gamma(1-\alpha)}\{\frac{x(t-\omega)}{\omega^\alpha}-\frac{x(t_0)}{(t-t_0)^\alpha}-\alpha\int_{t_0}^{t-\omega}(t-\tau)^{-\alpha-1}x(\tau)d\tau\}&
\text{$t>t_0+\omega$},
\end{cases}
\end{split}
\end{equation}
where $0<\alpha<1$ and $x(t_0)=x_0$. It is easy to find that the function
\begin{equation}\begin{split}\
g(x,t)=\begin{cases}
f(x,t)&\text{$t_0\leq t\leq t_0+\omega$},\\
f(x,t)+\frac{1}{\Gamma(1-\alpha)}\{\frac{x(t-\omega)}{\omega^\alpha}-\frac{x(t_0)}{(t-t_0)^\alpha}-\alpha\int_{t_0}^{t-\omega}(t-\tau)^{-\alpha-1}x(\tau)d\tau\}&
\text{$t>t_0+\omega$},
\end{cases}
\end{split}
\end{equation}
is continuous in $t$ and $x$ on $[t_0,+\infty)\times R^n$. Then according to \cite{2}, there exists a continuous solution $x(t)$ to the Cauchy type problem $(25)$ on the interval $[t_0,+\infty)$. Thus, there exists a continuous solution $x(t)$ to the Cauchy type problem $(24)$ on the interval $[t_0,+\infty)$. This means that the integrability of $f(x,t)$ and the continuity of operator ${_{s_\omega(t)}\tilde{D}^\alpha_t}$ can ensure the existence and continuity of solutions to the system $(24)$. But it should be pointed out that establishing the uniqueness conditions of solutions to the system $(24)$ is not the focus of this paper. Therefore, instead of discussing the uniqueness of solutions in this article, we assume that there exists a unique solution to the Cauchy type problem $(24)$.

Similar to $Definition~3$, the equilibrium point of $(24)$ is defined as follows.

\textbf{Definition 6 }: The $n$-dimensional constant vector $x^*$ is an equilibrium point of the short memory fractional system $(24)$, if $f(x^*,t)=0$.

Then Lyapunov direct method is extended to the case of short memory fractional systems and the asymptotic stability conditions are obtained.

\textbf{Theorem 4 }: Let $x^*=0$ be an equilibrium point for the non-autonomous system $(24)$. Assume that there exist
a Lyapunov function $V(t,x(t))$, a positive constant $\lambda$ and class-$K$ functions $\phi_i$ $(i=1,2)$ satisfying
\begin{equation}
\begin{aligned}
\phi_1(\parallel x\parallel)\leq V(t,x)\leq \phi_2(\parallel x\parallel),
\end{aligned}
\end{equation}
\begin{equation}
\begin{aligned}
{_{s_\omega(t)}\tilde{D}^\alpha_t}V(t,x(t))\leq-\lambda V(t,x(t)),
\end{aligned}
\end{equation}
and
\begin{equation}
\begin{aligned}
\lambda>\frac{1}{\omega^\alpha\Gamma(1-\alpha)},
\end{aligned}
\end{equation}
where $\Gamma(\cdot)$ is Gamma function. Then the equilibrium point $x^*=0$ of the system $(24)$ is Lyapunov asymptotically stable.

\textbf{Proof }: Consider the Caputo fractional differential equation with time delay
\begin{equation}
  \left\{
\begin{aligned}
&{^C_{t_0}D^\alpha_t}W(t)= -\lambda W(t)+\frac{1}{\omega^\alpha\Gamma(1-\alpha)}W(t-\omega),\\
&W(s)=V(t_0,x(t_0)),~~s\in[t_0-\omega,t_0],
  \end{aligned}
\right.
\end{equation}
where $W(t)\geq 0$ and $W(t)\in R$. According to $Lemma~5$, the zero equilibrium point of the system $(30)$ is Lyapunov asymptotically stable. And using $Theorem~3$, we have
\begin{equation}
\begin{aligned}
0\leq V(t,x)\leq W(t),~~t\geq t_0.
\end{aligned}
\end{equation}
Obviously, $\lim\limits_{t\rightarrow \infty}V(t,x)=0$ because $\lim\limits_{t\rightarrow \infty}W(t)=0$. Thus the equilibrium point $x^*=0$ of system $(24)$ is attractive.

According to the stability of the system $(30)$, we can obtain that, for any constant $\varepsilon_1>0$, there exists a positive constant $\delta_1$ such that $W(t)\leq \varepsilon_1$ when $W(t_0)\leq \delta_1$. Then for any constant $\varepsilon_2>0$, there exits a positive constant $\delta_2$ such that $V(t,x)\leq W(t)\leq \varepsilon_2$ when $V(t_0,x(t_0))\leq \delta_2$, where $\delta_2$ is an element in the set $\{\mu\mid W(t_0)\leq \mu, ~W(t)\leq\varepsilon_2,~t\geq t_0,~\mu>0\}$ and the initial condition $W(t_0)=V(t_0,x(t_0))$ in $(30)$ is used here. Hence the solution $V(t,x(t))$ of the differential inequality $(28)$ is stable. From the inequality $(27)$, the equilibrium point $x^*=0$ of system $(24)$ is stable. Therefore, the equilibrium point $x^*=0$ of system $(24)$ is Lyapunov asymptotically stable.

This completes the proof.

\textbf{Remark 6 }: In $Theorem~4$, the stability of the zero equilibrium point of $(24)$ is studied. Then the case of $x^*\neq 0$ can be discussed as follows. If $x^*\neq0$ is an equilibrium point for the non-autonomous system $(24)$, then the change of variables $\bar{x}=x-x^*$ can be applied and the following equality can be obtained,
\begin{equation}
\begin{aligned}
{_{s_\omega(t)}\tilde{D}^\alpha_t}\bar{x}(t)= f(\bar{x}+x^*,t)=g(\bar{x},t),
\end{aligned}
\end{equation}
where $g(0,t)=0$. Obviously, the system $(32)$ has an equilibrium at the origin and the stability analysis of nonzero equilibrium points of the system $(24)$ is equivalent to the stability analysis of zero equilibrium point of the system $(32)$. Therefore $Theorem~4$ is still valid when analyzing the stability of nonzero equilibrium points of the system $(24)$.

\textbf{Remark 7 }: It is easy to find that the larger the parameter $\omega$, the easier it is to satisfy the condition $(29)$. In particular, when $\omega\rightarrow+\infty$, the equation $(24)$ can be rewritten as
\begin{equation}
\begin{aligned}
{^C_{t_0}D^\alpha_t}x(t)=f(x,t),~~0<\alpha<1.
\end{aligned}
\end{equation}
According to $Lemma~6$, the system $(33)$ is Lyapunov asymptotically stable when the conditions $(27)$ and $(28)$ are satisfied.

\textbf{Remark 8 }: Note that the system $(24)$ is local Lyapunov asymptotic stable when the conditions $(27)$, $(28)$ and $(29)$ are satisfied. If an additional condition that $\lim\limits_{\gamma\rightarrow \infty}\phi_1(\gamma)=\infty$ is met, then the system $(24)$ is global Lyapunov asymptotic stable.

Next, a theorem is proposed, by which the stability of some special differential systems can be verified directly.

\textbf{Theorem 5 }: Let $x(t)=(x_1(t),x_2(t),\cdots,x_n(t))^T\in R^n$ and $x^*=0$ be an equilibrium point for the non-autonomous system $(24)$. Assume that there exist a constant $\varphi>0$ and positive integers $m_i$ $(i=1,2,\cdots,n)$ satisfying
\begin{equation}
\begin{aligned}
\zeta^Tf(x,t)\leq-\varphi \sum\limits_{i=1}^nx_i^{2^{m_{i}}}
\end{aligned}
\end{equation}
and
\begin{equation}
\begin{aligned}
\varphi 2^{\hat{m}}>\frac{1}{\omega^\alpha\Gamma(1-\alpha)},
\end{aligned}
\end{equation}
where $\zeta=(x_1^{2^{m_1}-1}(t),x_2^{2^{m_2}-1}(t),\cdots,x_n^{2^{m_n}-1}(t))^T$ and $\hat{m}=\min\{m_1,m_2,\cdots,m_n\}$. Then the equilibrium point $x^*=0$ of the system $(24)$ is Lyapunov asymptotically stable.

\textbf{Proof }: Consider a Lyapunov function as follows,
\begin{equation}
\begin{aligned}
V(t,x)=\sum\limits_{i=1}^n\frac{1}{2^{m_i}}x_i^{2^{m_i}}(t).
\end{aligned}
\end{equation}

By applying $Theorem~2$, the following inequality holds,
\begin{equation}
\begin{aligned}
{_{s_\omega(t)}\tilde{D}^\alpha_t}V(t,x)&=\sum\limits_{i=1}^n\frac{1}{2^{m_i}}{_{s_\omega(t)}\tilde{D}^\alpha_t}x_i^{2^{m_i}}(t)\\
&\leq \sum\limits_{i=1}^n\frac{1}{2^{m_i-1}}x_i^{2^{m_i-1}}(t){_{s_\omega(t)}\tilde{D}^\alpha_t}x_i^{2^{m_i-1}}(t)\\
&\leq \sum\limits_{i=1}^n\frac{1}{2^{m_i-2}}x_i^{2^{m_i-1}+2^{m_i-2}}(t){_{s_\omega(t)}\tilde{D}^\alpha_t}x_i^{2^{m_i-2}}(t)\\
&\cdots\\
&\leq \sum\limits_{i=1}^nx_i^{2^{m_i}-1}(t){_{s_\omega(t)}\tilde{D}^\alpha_t}x_i(t)\\
&=\zeta^Tf(x,t)\\
&\leq -\varphi \sum\limits_{i=1}^nx_i^{2^{m_{i}}}(t)\\
&=-\varphi 2^{\hat{m}}\sum\limits_{i=1}^n(\frac{1}{2^{\hat{m}}}x_i^{2^{m_{i}}}(t))\\
&\leq-\varphi 2^{\hat{m}} V(t,x).
\end{aligned}
\end{equation}

According $Theorem~4$, the equilibrium point $x^*=0$ of system $(24)$ is Lyapunov asymptotically stable.

This completes the proof.

\textbf{Remark 9 }: In particular, when $\omega\rightarrow+\infty$, the method proposed in $Theorem~5$ is still valid and the corresponding discussion has been given in \cite{15}.

\section{Numerical examples}
In this section, three examples are shown to verify the correctness of the obtained theoretical results.

\textbf{Example 1 }: A short memory fractional system is given as follows,
\begin{equation}
\begin{aligned}
{_{s_\omega(t)}\tilde{D}^\alpha_t}y(t)=-y(t),
\end{aligned}
\end{equation}
where $\alpha=0.95$, $\omega=5$, $y(t_0)=3$ and $t_0=0$. Then consider a Caputo fractional differential equation with time delay
\begin{equation}
  \left\{
\begin{aligned}
&{^C_{t_0}D^\alpha_t}x(t)= -x(t)+\frac{1}{\omega^\alpha\Gamma(1-\alpha)}x(t-\omega),~~0<\alpha<1,\\
&x(s)=3,~~s\in[t_0-\omega,t_0],
  \end{aligned}
\right.
\end{equation}
where $\alpha=0.95$, $\omega=5$ and $t_0=0$. According to $Theorem~3$, we can obtain that $x(t)\geq y(t)$ for $t\geq 0$, and this conclusion can be verified by $Figure~1$ which shows the evolution of the states of systems $(38)$ and $(39)$.

\textbf{Example 2 }: The following equation is considered,
\begin{equation}
  \left\{
\begin{aligned}
&{_{s_\omega(t)}\tilde{D}^\alpha_t}x_1(t)=x_2(t),\\
&{_{s_\omega(t)}\tilde{D}^\alpha_t}x_2(t)=-x_1(t)-x_2(t),
  \end{aligned}
\right.
\end{equation}
where $\alpha=0.95$, $\omega=5$, $t_0=0$ and the initial conditions are chosen as $x_1(t_0)=7, x_2(t_0)=-3$. Then a Lyapunov function is considered as follows,
\begin{equation}
\begin{aligned}
V(t)=3x_1^2(t)+2x_1(t)x_2(t)+2x_2^2(t).
\end{aligned}
\end{equation}
By applying $Theorem~2$, we have
\begin{equation}
\begin{aligned}
{_{s_\omega(t)}\tilde{D}^\alpha_t}V(t)&\leq4x_1(t){_{s_\omega(t)}\tilde{D}^\alpha_t}x_1(t)+2(x_1(t)+x_2(t))({_{s_\omega(t)}\tilde{D}^\alpha_t}(x_1(t)+x_2(t)))+2x_2(t){_{s_\omega(t)}\tilde{D}^\alpha_t}x_2(t)\\
&=-2(x_1^2(t)+x_2^2(t))\\
&=-\frac{1}{2}(4x_1^2(t)+4x_2^2(t))\\
&\leq- \frac{1}{2}(3x_1^2(t)+2x_1(t)x_2(t)+2x_2^2(t))\\
&=-\frac{1}{2}V(t)<0.
\end{aligned}
\end{equation}
From the inequality $(42)$, the short memory fractional derivative of the Lyapunov function $(41)$ is negative definite and the conditions $(28)$ and $(29)$ can be satisfied, thus it can be concluded from $Theorem~4$ that the origin of the system $(40)$ is asymptotically stable. $Figure~2$ shows the evolution of the states of the system $(40)$.

\textbf{Example 3 }: Consider the following differential equation,
\begin{equation}
  \left\{
\begin{aligned}
&{_{s_\omega(t)}\tilde{D}^\alpha_t}x_1(t)=-x_1(t)+x_2^3(t),\\
&{_{s_\omega(t)}\tilde{D}^\alpha_t}x_2(t)=-x_1(t)-x_2(t),
\end{aligned}
\right.
\end{equation}
where $\alpha=0.95$, $\omega=5$, $t_0=0$ and the initial conditions are chosen as $x_1(t_0)=3, x_2(t_0)=-5$. It is easy to find that
\begin{equation}
\begin{aligned}
\zeta^Tf(x_1,x_2)&=-x_1^2(t)-x_2^4(t),
\end{aligned}
\end{equation}
where $\zeta=(x_1(t),x_2^3(t))^T$ and $f(x_1,x_2)=(-x_1(t)+x_2^3(t),-x_1(t)-x_2(t))^T$. Then according to $Theorem~5$, the origin of the equation $(43)$ is asymptotically stable. And $Figure~3$ shows the evolution of the states of the system $(43)$.

\textbf{Remark 10 }: According to $Remark~5$, the solution of the short memory fractional equation $(24)$ is also a solution of the Caputo fractional system $(25)$. Hence, the numerical solution of the system $(24)$ can be obtained by using a predictor-corrector method \cite{19} to solve the numerical solution of the equation $(25)$. And the numerical solutions of equations $(38)(40)(43)$ are obtained in the same way.

\begin{figure}
    \begin{center}
        \includegraphics[width=0.95\linewidth]{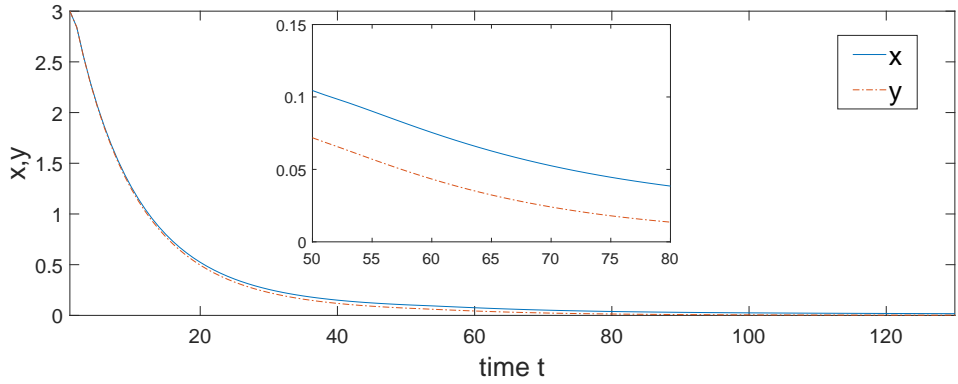}
       \caption{Evolution of the states of systems $(38)$ and $(39)$. }\label{Figure:1}
    \end{center}
    \end{figure}

    \begin{figure}
    \begin{center}
        \includegraphics[width=0.95\linewidth]{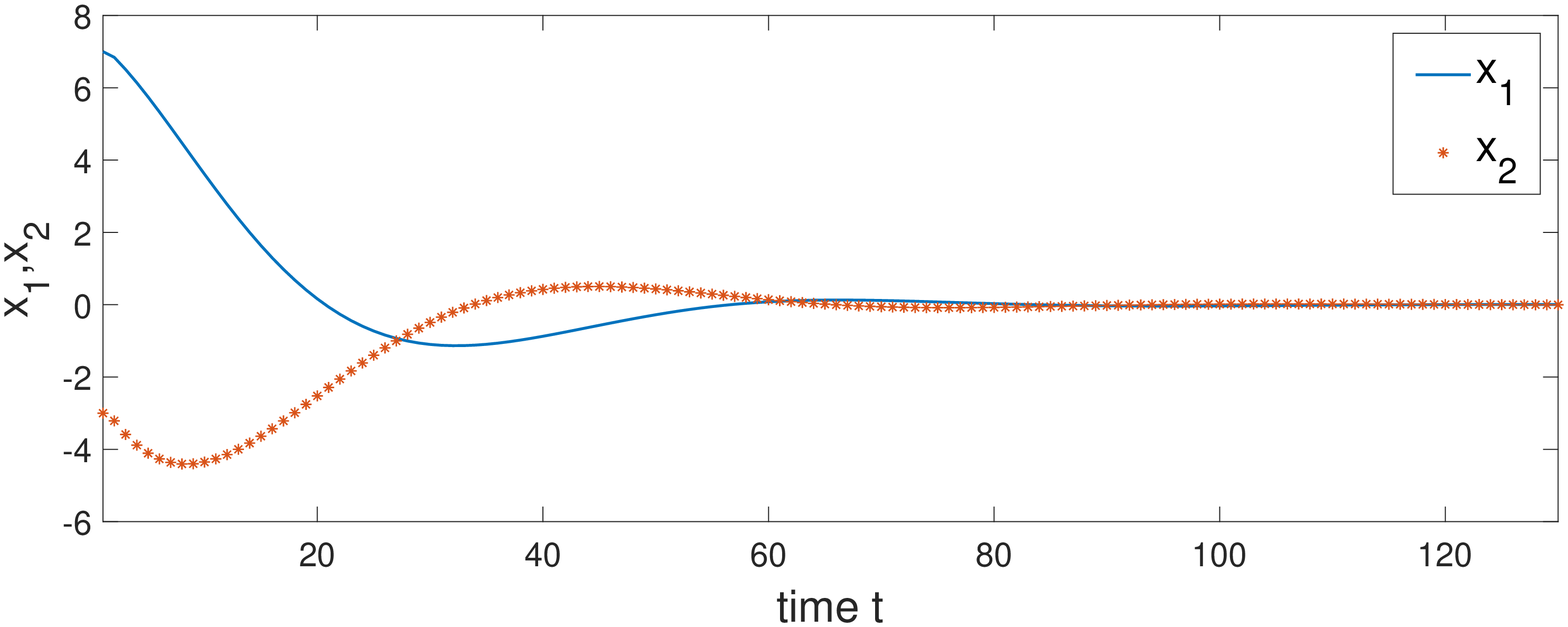}
       \caption{Evolution of the states $x_1(t)$ and $x_2(t)$ of the system $(40)$. }\label{Figure:2}
    \end{center}
    \end{figure}

\begin{figure}
    \begin{center}
        \includegraphics[width=0.95\linewidth]{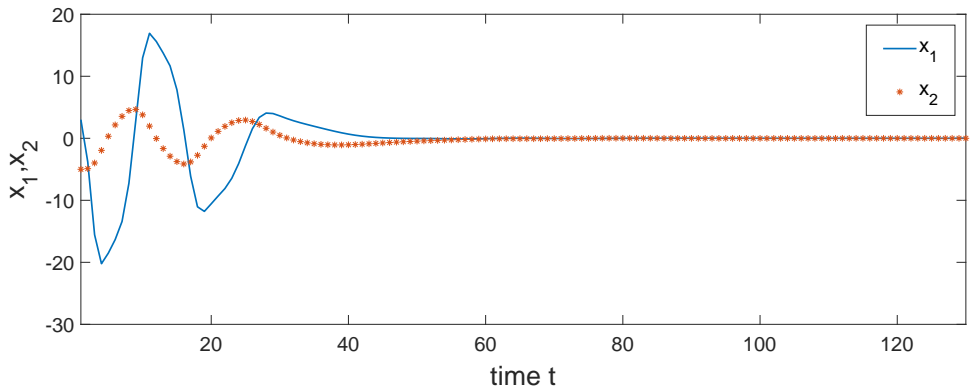}
       \caption{Evolution of the states $x_1(t)$ and $x_2(t)$ of the system $(43)$. }\label{Figure:2}
    \end{center}
    \end{figure}

\section{Conclusion}
In this manuscript, the stability analysis of short memory fractional systems is investigated. Firstly, a new short memory fractional differential operator is defined and some properties of this operator are derived. Then a comparison theorem for short memory fractional equations is proposed, based on which some stability conditions of short memory fractional systems are obtained. In addition, a special theorem is presented, which can make it easier to judge the stability of some special systems. Finally, several numerical examples are given to verify the effectiveness of the main results.

And the stability conditions of short memory fractional systems with time delays will be studied in our next work.

\section*{Acknowledgment}
This work is supported by the Natural Science Foundation of Beijing Municipality (No. Z180005), the National Nature Science Foundation of China (No. 61772063) and the National Nature Science Foundation of China (No. 61973329).

\end{document}